\documentclass[12pt]{amsart}
\usepackage{amssymb}

\long\def\forget#1\forgotten{}
\newcommand{\issuenumber}{20}
\newcommand{\issuemonth}{March}
\newcommand{\issueyear}{2007}

\newtheorem{thm}{Theorem}[section]
\newtheorem{prob}[thm]{Problem}

\newtheorem{issue}{Issue}

\theoremstyle{definition}

\theoremstyle{remark}

\newcommand{\ed}{
\newpage

\section{Unsolved problems from earlier issues}

\begin{issue}
Is $\binom{\Omega}{\Gamma}=\binom{\Omega}{\Tau}$?
\end{issue}

\begin{issue}
Is $\ufin(\cO,\Omega)=\sfin(\Gamma,\Omega)$?
And if not, does $\ufin(\cO,\Gamma)$ imply
$\sfin(\Gamma,\Omega)$?
\end{issue}

\stepcounter{issue}

\begin{issue}
Does $\sone(\Omega,\Tau)$ imply $\ufin(\Gamma,\Gamma)$?
\end{issue}

\begin{issue}
Is $\fp=\fp^*$? (See the definition of $\fp^*$ in that issue.)
\end{issue}

\begin{issue}
Does there exist (in ZFC) an uncountable set satisfying $\sone(\BG,\B)$?
\end{issue}

\stepcounter{issue}

\begin{issue}
Does $X \nin \NON(\M)$ and $Y\nin\mathsf{D}$ imply that
$X\cup Y\nin \COF(\M)$?
\end{issue}

\begin{issue}[CH]
Is $\split(\Lambda,\Lambda)$ preserved under finite unions?
\end{issue}

\begin{issue}
Is $\cov(\M)=\fo$? (See the definition of $\fo$ in that issue.)
\end{issue}

\begin{issue}
Does $\sone(\Gamma,\Gamma)$ always contain an element of cardinality $\fb$?
\end{issue}

\begin{issue}
Could there be a Baire metric space $M$ of weight $\aleph_1$ and a partition
$\mathcal{U}$ of $M$ into $\aleph_1$ meager sets where for each ${\mathcal U}'\subset\mathcal U$,
$\bigcup {\mathcal U}'$ has the Baire property in $M$?
\end{issue}

\stepcounter{issue} 

\begin{issue}
Does there exist (in ZFC) a set of reals $X$ of cardinality $\fd$ such that all
finite powers of $X$ have Menger's property $\sfin(\cO,\cO)$?
\end{issue}

\begin{issue}
Can a Borel non-$\sigma$-compact group be generated by a Hurewicz subspace?
\end{issue}

\begin{issue}[MA]
Is there an uncountable $X\sbst\R$ satisfying $\sone(\BO,\BG)$?
\end{issue}

\begin{issue}[CH]
Is there a totally imperfect $X$ satisfying $\ufin(\cO,\Gamma)$
that can be mapped continuously onto $\Cantor$?
\end{issue}

\begin{issue}[CH]
Is there a Hurewicz $X$ such that $X^2$ is Menger but not Hurewicz?
\end{issue}

\begin{issue}
Does the Pytkeev property of $C_p(X)$ imply the Menger property of $X$?
\end{issue}

\begin{issue}
Does every hereditarily Hurewicz space satisfy $\sone(\BG,\BG)$?
\end{issue}

\general\end{document}}

\newcommand{\Cantor}{{\{0,1\}^\N}}

\newcommand{\fb}{\mathfrak{b}}

\newcommand{\fd}{\mathfrak{d}}
\newcommand{\fp}{\mathfrak{p}}

\newcommand{\NON}{{\mathsf   {NON}}}
\newcommand{\COF}{{\mathsf   {COF}}}

\newcommand{\M}{\mathcal{M}}

\newcommand{\cov}{\mathsf{cov}}

\newcommand{\R}{\mathbb{R}}
\newcommand{\Q}{\mathbb{Q}}

\newcommand{\fo}{\mathfrak{od}}

\renewcommand{\split}{\mathsf{Split}}
\newcommand{\bq}{\begin{quote}}
\newcommand{\eq}{\end{quote}}
\newcommand{\cO}{\mathcal{O}}
\newcommand{\B}{\mathcal{B}}
\newcommand{\BG}{\B_\Gamma}

\newcommand{\BO}{\B_\Omega}
\newcommand{\CG}{C_\Gamma}

\newcommand{\sone}{\mathsf{S}_1}    \newcommand{\sfin}{\mathsf{S}_{fin}}

\newcommand{\ufin}{\mathsf{U}_{fin}}

\newcommand{\nin}{\not\in}

\newcommand{\NN}{{\N^\N}}
\newcommand{\N}{\mathbb{N}}

\newcommand{\sbst}{\subseteq}
\newcommand{\by}[2]{\par\hfill\emph{#1}, #2}
\newcommand{\nby}[1]{\par\hfill\emph{#1}}
\newcommand{\Tau}{\mathrm{T}}
\newcommand{\CE}{\textsc{CE}}

\newcommand{\be}{\begin{enumerate}}
\newcommand{\ee}{\end{enumerate}}
\newcommand{\bi}{\begin{itemize}}
\newcommand{\ei}{\end{itemize}}
\newcommand{\itm}{\item}

\newcommand{\general}{\small\vfill\par\noindent\hrulefill\par
\noindent\textbf{Previous issues.} The previous issues of this
bulletin, which contain general information (first issue), basic
definitions, research announcements, and open problems (all
issues) are available online,
at \texttt{http://front.math.ucdavis.edu/search?\&t=\%22SPM+Bulletin\%22}
\\[0.1cm]
\textbf{Contributions.}
Please submit your contributions (announcements, discussions, and open problems)
by e-mailing us. It is preferred to write them
in \LaTeX{}.
The authors are urged to use as standard notation as possible, or otherwise give
the definitions or a reference to where the notation is explained.
Contributions to this bulletin would not require any transfer of copyright,
and material presented here can be published elsewhere.\\[0.1cm]
\textbf{Subscription.}
To receive this bulletin (free) to your
e-mailbox, e-mail us:\\
{boaz.tsaban@weizmann.ac.il}
}

\newcommand{\nArxPaper}[5]{\subsection{#2}{#4}\par\hfill{\arx{#1}}\par\hfill\emph{#3}}

\newcommand{\nAMSPaper}[4]{\subsection{#2}{#4}\par\hfill{\texttt{#1}}\par\hfill\emph{#3}}

\newcommand{\arx}[1]{\texttt{http://arxiv.org/math/#1}}
\newcommand{\url}[1]{\bq\texttt{#1}\eq}
\newcommand{\online}[1]{The paper is available online at \url{#1}}

\title[$\mathcal{SPM}$ Bulletin \textbf{\issuenumber} (\issuemonth{} \issueyear)]{%
$\mathcal{SPM}$ Bulletin\\[0.5cm]
Issue number \issuenumber: \issuemonth{} \issueyear{} \CE{}}

\begin{document}
\maketitle

\tableofcontents

\section{Editor's note}

The \emph{Third Workshop on Coverings, Selections, and Games in Topology}
is just around the corner. Visit
\url{http://www.pmf.ni.ac.yu/spm2007/index.html}
for details. We are looking forward to this opportunity to meet
old and new friends and colleagues, and to exchange some new results and
ideas.

\medskip

\by{Boaz Tsaban}{boaz.tsaban@weizmann.ac.il}

\hfill \texttt{http://www.cs.biu.ac.il/\~{}tsaban}

\section{Research announcements}

\nArxPaper{math.GN/0701249}
{On the Pytkeev property in spaces of continuous functions (II)}
{Boaz Tsaban and Lyubomyr Zdomskyy}
{We prove that for each Polish space $X$, the space $C(X)$ of continuous
real-valued functions on $X$ satisfies a strong version of the Pytkeev property,
if endowed with the compact-open topology. (This shows that whereas it need not
be metrizable, it is ``very close'' to that.)

We also consider the Pytkeev property in the case where $C(X)$ is endowed with
the topology of pointwise convergence.}

\subsection{Partial order embeddings with convex range}
{
A careful study is made of embeddings of posets which have a convex range.
We observe that such embeddings
share nice properties
with the homomorphisms of more restrictive categories; for example, we show that
every order embedding between two lattices with convex range
is a continuous lattice homomorphism.
A number of posets are considered; for example, we prove that
every product order embedding $\sigma:\NN\to\NN$ with convex range is of the
form
\begin{equation}
  \sigma(x)(n)=\bigl((x\circ g_\sigma)+y_\sigma\bigr)(n)
  \quad\text{if $n\in K_\sigma$},\label{eq:68}
\end{equation}
and $\sigma(x)(n)=y_\sigma(n)$ otherwise, for all $x\in\NN$, where
$K_\sigma\subseteq\N$, $g_\sigma:K_\sigma\to\N$ is a bijection and
$y_\sigma\in\NN$. The most complex poset examined here is the quotient of the
lattice of Baire measurable functions,
with codomain of the form $\N^I$ for some index set $I$, modulo equality on a
comeager subset of the domain, with its `natural' ordering.
}

\noindent\texttt{http://homepage.univie.ac.at/James.Hirschorn/research/embeddings/}\\
\texttt{embeddings.html}
\nby{James Hirschorn}

\nArxPaper{math.GN/0702296}
{Resolvability vs.\ almost resolvability}
{Istvan Juhasz and Saharon Shelah and Lajos Soukup}
{
 A space $X$ is $\kappa$-resolvable (resp.\ almost $\kappa$-resolvable) if it contains
$\kappa$ dense sets that are pairwise disjoint (resp. almost disjoint over the
ideal of nowhere dense subsets of $X$). Answering a problem raised by Juhasz,
Soukup, and Szentmiklossy, and improving a consistency result of Comfort and
Hu, we prove, in ZFC, that for every infinite cardinal ${\kappa}$ there is an
almost $2^{\kappa}$-resolvable but not $\aleph_1$-resolvable space of dispersion
character ${\kappa}$.
}

\nArxPaper{math.LO/0702295}
{Lindel\"of spaces of singular density}
{Istvan Juhasz and Saharon Shelah}
{
  A cardinal $\lambda$ is called omega-inaccessible if for all $\mu < \lambda$ we have
$\mu^\omega<\lambda$. We show that for every $\omega$-inaccessible cardinal $\lambda$
there is a CCC (hence cardinality and cofinality preserving) forcing that adds
a hereditarily Lindel\"of regular space of density $\lambda$. This extends an
analogous earlier result of ours that only worked for regular $\lambda$.
}

\nArxPaper{math.LO/0702309}
{$P(\omega)/{\rm fin}$ and projections in the Calkin algebra}
{Eric Wofsey}
{
  We investigate the set-theoretic properties of the lattice of projections in
the Calkin algebra of a separable infinite-dimensional Hilbert space in
relation to those of the Boolean algebra $P(\omega)/{\rm fin}$, which is
isomorphic to the sublattice of diagonal projections. In particular, we prove
some basic consistency results about the possible cofinalities of well-ordered
sequences of projections and the possible cardinalities of sets of mutually
orthogonal projections that are analogous to well-known results about
$P(\omega)/{\rm fin}$.
}

\subsection{Everywhere meagre and everywhere null sets}
We introduce new classes of small subsets of the reals, having natural combinatorial definitions, namely
everywhere meagre and everywhere null sets, which lie between the $\sigma$-ideal ${\mathcal I}_0$, introduced by
Repick\'y in [M. Repick\'y,
\emph{Mycielski ideal and the perfect set theorem},
Proc.\ AMS \textbf{132} (2004), 2141--2150]
and the $\sigma$-ideal ${\mathbb B}_2$, introduced by Ros\l{}anowski in [A. Ros\l{}anowski,
\emph{On game ideals},
Colloq.\ Math.\ \textbf{59} (1990), 159--168.].
We investigate properties of these sets, in particular we
show that these classes are closed under taking products and projections. We also prove several relations
between these classes and other well-known classes of small subsets of the reals.

\hfill\texttt{http://www.math.uni.wroc.pl/\~{}kraszew/files/papers.html}

\nby{Jan Kraszewski}

\nAMSPaper{http://www.ams.org/journal-getitem?pii=S0002-9939-07-08714-X}
{Transversals for strongly almost disjoint families}
{Paul J. Szeptycki}
{For a family of sets $A$, and a set $X$, $X$ is said to be a transversal of $A$ if
$X\subseteq \bigcup A$ and $|a\cap X|=1$ for each $a\in A$. $X$ is said to be a
Bernstein set for $A$ if $\emptyset\not=a\cap X\not=a$ for each $a\in A$. When an
almost disjoint family admits a set like a transversal or Bernstein set was first
studied in [P.\ Erd\"os and A.\ Hajnal \emph{On a property of families of sets},
Acta Math.\ Acad.\ Sci.\ Hung., 12 (1961) 87--124.].
In this note we introduce the following notion: a family of sets
$A$ is said to admit a $\sigma$-transversal if $A$ can be written as $A=\bigcup\{A_n:n\in \omega\}$
such that each $A_n$ admits a transversal. We study the question when an almost disjoint family
admits a $\sigma$-transversals and related questions.}

\nArxPaper{math.LO/0703091}
{A family of covering properties for forcing axioms and strongly compact cardinals}
{Matteo Viale}
{This paper presents the main results in my Ph.D. thesis.
Several proofs of SCH are presented introducing a family of covering properties
which implies both SCH and the failure of various forms of square. These
covering properties are also applied to investigate models of strongly compact
cardinals or of strong forcing axioms like MM or PFA.}

\subsection{Splitting families and forcing}

According to [M.S. Kurili\'c, \emph{Cohen-stable families of subsets of the integers},
J. Symbolic Logic \textbf{66} (2001) 257--270], adding
a Cohen real destroys a splitting family $S$ on $\N$ if and only if S is isomorphic to a splitting family on
the set of rationals, $\Q$,
whose elements have nowhere dense boundaries.
Consequently, $|S| < \cov(\M)$ implies the Cohen-indestructibility of $S$. Using
the methods developed in [J. Brendle, S. Yatabe, \emph{Forcing indestructibility of MAD families},
Ann.\ Pure Appl.\ Logic \textbf{132} (2005) 271--312],
the stability of splitting families in several forcing extensions is characterized in a similar way
(roughly speaking, destructible families have members with ``small generalized boundaries'' in the
space of the reals). Also, it is proved that a splitting
family is preserved by the Sacks (respectively: Miller, Laver) forcing if and only if it is preserved by
some forcing which adds a new (respectively: an unbounded, a dominating) real.
The corresponding hierarchy of splitting families is investigated.

\hfill\texttt{http://dx.doi.org/10.1016/j.apal.2006.08.002}

\nby{Milo\v{s} S. Kurili\'c}

\nArxPaper{math.LO/0703302}
{A Sacks Real out of Nowhere}
{Jakob Kellner, Saharon Shelah}
{There is a proper countable support iteration of length $\omega$ adding no
new reals at finite stages and adding a Sacks real in the limit.}

\nArxPaper{math.GR/0703304}
{Reflection principle characterizing groups in which unconditionally closed sets are algebraic}
{Dikran Dikranjan, Dmitri Shakhmatov}
{We give a necessary and sufficient condition, in terms of a certain
reflection principle, for every unconditionally closed subset of a group $G$ to
be algebraic. As a corollary, we prove that this is always the case when $G$ is a
direct product of an Abelian group with a direct product (sometimes also called
a direct sum) of a family of countable groups. This is the widest class of
groups known to date where the answer to the 63 years old problem of Markov
turns out to be positive. We also prove that whether every unconditionally
closed subset of $G$ is algebraic or not is completely determined by countable
subgroups of $G$.
}

\nArxPaper{math.GR/0703397}
{Unconditionally $\tau$-closed and $\tau$-algebraic sets in groups}
{Ol'ga V. Sipacheva}
{Families of unconditionally $\tau$-closed and $\tau$-algebraic sets in a group are
defined. It is proved that, if any element of a group $G$ has at most $\tau$
conjugates, then these families coincide in $G$. It follows that any
unconditionally closed set is algebraic in a group in which every element has
at most countably many conjugates.}

\subsection{Stratifiability of $C_k(X)$ for a class of separable metrizable $X$}
Let $X$ be a separable metrizable space.
It is proved that the space
$C_k(X)$
of all continuous real-valued functions on $X$ with the compact-open topology
is stratifiable if and only if $X$ is Polish.

\hfill\texttt{http://erezn.andante.ru/topology/sckx3.ps}

\nby{E. A. Reznichenko}

\nArxPaper{math.GN/0703429}
{Dissipated Compacta}
{Kenneth Kunen}
{The dissipated spaces form a class of compacta which contains both the
scattered compacta and the compact LOTSes (linearly ordered topological
spaces), and a number of theorems true for these latter two classes are true
more generally for the dissipated spaces. For example, every regular Borel
measure on a dissipated space is separable.

  A product of two compact LOTSes is usually not dissipated, but it may satisfy
a weakening of that property. In fact, the degree of dissipation of a space can
be used to distinguish topologically a product of $n$ LOTSes from a product of $m$
LOTSes.
}

\nArxPaper{math.LO/0703647}
{A note on strong negative partition relations}
{Todd Eisworth}
{We analyze a natural function definable from a scale at a singular cardinal,
and using this function we are able to obtain quite strong negative
square-brackets partition relations at successors of singular cardinals. The
proof of our main result makes use of club-guessing, and as a corollary we
obtain a fairly easy proof of a difficult result of Shelah connecting weak
saturation of a certain club-guessing ideal with strong failures of
square-brackets partition relations. We then investigate the strength of weak
saturation of such ideals and obtain some results on stationary reflection.
}

\nArxPaper{math.GN/0703728}
{First countable spaces without point-countable $\pi$-base}
{Istvan Juhasz, Lajos Soukup, and Zoltan Szentmiklossy}
{We answer several questions of V. Tka\v{c}uk from [Point-countable
$\pi$-bases in first countable and similar spaces, Fund. Math. 186 (2005), pp.
55--69.] by showing that
\be
\itm There is a ZFC example of a first countable, 0-dimensional Hausdorff
space with no point-countable $\pi$-base (in fact, the order of any $\pi$-base
of the space is at least $\aleph_\omega$);
\itm If there is a $\kappa$-Suslin line then there is a first countable GO
space of cardinality $\kappa^+$ in which the order of any $\pi$-base is at
least $\kappa$;
\itm It is consistent to have a first countable, hereditarily Lindel\" of
regular space having uncountable $\pi$-weight and $\omega_1$ as a caliber (of
course, such a space cannot have a point-countable $\pi$-base).
\ee
}

\nArxPaper{math.GR/0703726}
{Less than continuum many translates of a compact nullset may cover any infinite profinite group}
{Miklos Abert}
{We show that it is consistent with the axioms of set theory that every
infinite profinite group G possesses a closed subset X of Haar measure zero
such that less than continuum many translates of X cover G. This answers a
question of Elekes and Toth and by their work settles the problem for all
infinite compact topological groups.}

\nArxPaper{math.GN/0703835}
{Projective $\pi$-character bounds the order of a $\pi$-base}
{Istvan Juhasz and Zoltan Szentmiklossy}
{All spaces below are Tychonov. We define the projective $\pi$-character $p(X)$ of
a space $X$ as the supremum of the values $\pi\chi(Y)$ where $Y$ ranges over all
continuous images of $X$. Our main result says that every space $X$ has a $\pi$-base
whose order is at most $p(X)$, that is every point in $X$ is contained in at most
$p(X)$-many members of the $\pi$-base. Since $p(X)$ is at most $t(X)$ for compact $X$,
this provides a significant generalization of a celebrated result of
Shapirovskii.}

\section{Problem of the Issue}

We recall from \cite{BRR1} that a topological space $X$ is  a
\emph{$QN$-space} (resp \emph{wQN-space}), if  every pointwise
convergent to $0$ sequence $(f_n:X\to\R)_{n\in\N}$ of continuous
functions (contains a subsequence which) converges quasi-normally,
which means that there exists a convergent to $0$ sequence
$(\epsilon_n)_{n\in\N}$ of positive reals, such that for every $x\in X$
the inequality $|f_n(x)|<\epsilon_n$ holds for all but finitely many
$n\in\N$. Of course, every QN-space is a wQN-space. In
addition, every wQN-space has the Hurewicz property
$\ufin(\cO,\Gamma)$, see
\cite[Corollary~2,2]{BRR1}. It is also known \cite{Sa07, BH07} that $X$
is a wQN-space if, and only if, it satisfies the selection hypothesis
$\sone(\CG, \CG)$, and being a QN-spaces is  equivalent
\cite{BH07} to a certain selection principle stronger than
$\sone(\Gamma,\Gamma)$. Combining these results we conclude that the
above properties are related as follows:
$${QN} \Rightarrow
\sone(\Gamma,\Gamma) \Rightarrow \sone(\CG,\CG)
\Rightarrow \ufin(\cO,\Gamma)$$
Neither the first nor the third implication can be reversed in ZFC,
even for sets of reals without perfect subsets, see \cite{BRR1}
and \cite{TsZdPAMS} respectively. The question whether the  second
implication can be reversed was posed in \cite{ScEWJM} and is
still open.

 Passing to hereditary versions of these properties, we have
the following implications:
$$
QN
\stackrel{\mbox{\cite{Re}}}{=}
hQN
\Rightarrow
h\sone(\Gamma,\Gamma)
\stackrel{\mbox{\cite{Ha}}}{=}
h\sone(\CG,\CG)
\Rightarrow
h\ufin(\cO,\Gamma).$$
It is open whether any of the two displayed implications
can be reversed.
The question whether every hereditarily Hurewicz space a QN-space
is particularly interesting.

As it was recently showed in \cite{TsZdnew}, a space $X$ is a
QN-space if, and only if, it has the property $\sone(\BG,\BG)$,
which is equivalent \cite{TsSc} to the Hurewicz
property applied to the family of all countable Borel covers of
$X$, and to the property that all Borel images of $X$ in $\NN$ are bounded \cite{TsSc}.
In light of this equivalence the above question can be
reformulated as follows:
\begin{prob}
Does every hereditarily Hurewicz space satisfy $\sone(\BG,\BG)$?
\end{prob}

\nby{Lyubomyr Zdomskyy}

\ed